\documentclass{article}
\usepackage{amssymb, amsmath}
\begin{document}
\title{ A Quick Overview On Regular Algebraic $K$-Theory For Groups }
\author{U. Haag}
\date{Feb. 2014\\ \texttt{ Contact: haag@mathematik.hu-berlin.de}}
\maketitle
\large{
\noindent
The purpose of this expos\'{e} is to informally introduce the reader of the
manuscript  \emph{"Regular Algebraic $K$-Theory For Groups"} which develops the theory from scratch to its basic ideas and motivations. $K$-theory in all of its variants, topological, algebraic and 
operator theoretical, has seen a rapid development from the beginnings, starting with the work of [Grothendieck], [Atiyah] and others, which evolutionary process is still continuing to this day. It is now a widely accepted practice that a mathematician of some renomm\'{e}e should have his own Algebraic 
$K$-theory. Besides the most common definition of the higher algebraic $K$-groups due to [Quillen] there are theories of [Milnor], [Swan], [Karoubi-Villamayor] and others whose not being mentioned explicitely should be attributed to the authors ignorance. 
One should not forget to mention the overall important contribution of [Schur]. Indeed, Regular Algebraic $K$-theory may be properly seen as a theory of higher Schur multipliers.
Recently also bivariant Algebraic 
$K$-theories have been constructed along the lines and in the spirit of the ingenious conception of Operator $KK$-theory by [Kasparov], one for topological (locally convex) algebras due to [Cuntz], and a version which applies in a completely algebraic setting by [Cortinas-Thom]. 
One might ask why there should be any desire for yet another variant of Algebraic $K$-theory. I will try to give some motivation for 
the development of this theory which is still incomplete in many ways. First of all the most "basic" 
$K$-groups $\, K_0 ( R )\, ,\, K_1 ( R )\, $ of a (discrete) ring $\, R\, $ (and also $\, K_2 ( R )\, $ by Milnor's construction) can be defined in purely algebraic terms without referring to such concepts as homotopy groups or spectra. Of course it is commonly seen as an advantage to have such concepts at ones disposal which facilitate the computation of certain groups and are useful in many abstract considerations. It seems however that a theory which is encoded in purely algebraic terms might turn out to be more "universal" in the sense that it is more easily and naturally transformed to some related functor and the essentially algebraic dependence of the functor on the (algebraic) structure of the object under consideration appears with greater clarity. Another interesting aspect of Regular Algebraic $K$-theory is that it relies on the linking of the category of (discrete) rings with the (somewhat simpler or less structured) category of (discrete) groups. 
Namely it uses the natural transformation from rings to groups which assigns to a given unital ring 
$\, R\, $ its associated group generated by infinitedimensional elementary matrices $\, E ( R )\, $ with entries in $\, R\, $. Thus the theory is really defined on the category of groups rather than rings but by the substitution 
$\, R \rightsquigarrow E ( R )\, $ it can be applied to the category of rings. This gives the theory a certain flexibility and generality since it may also be applied in many other contexts where groups occur associated with some geometric or algebraic structure. This functor from rings to groups is already inherent in Quillen's 
$+$-construction and it can be shown that Quillen's higher $K$-groups depend only on $\, E ( R )\, $, so that for the purpose of comparison with Quillen's theory nothing is given away through this passage.
Of course there is also a natural transformation from groups to rings by assigning to a given group $\, G\, $ its group ring $\, CG\, $ over some chosen coefficient ring 
$\, C\, $ which makes possible to define a homology theory for groups on considering the Algebraic 
$K$-groups (of some sort) of the associated group ring. Using our definition of $K$-groups this amounts to a "stabilized" version of the theory on considering $\, K_* ( E ( CG ) )\, $ (or $\, K_* ( GL ( CG ) )\, $) instead of $\, K_* ( G )\, $. While this stabilization might be of some interest in certain situations it appears that one looses a great deal of generality and flexibility if starting with such a definition from the beginning.
We proceed by describing the most basic features of the construction. 
There are two variants of the theory, the first is named \emph{$K^J$-theory} and the second is \emph{Regular $K$-theory}. The former is the more general of the two since it is defined on the category of \emph{all} normal subgroups, i.e. on the category of \emph{normal pairs} $(N, G)$ where $\, N\, $ is a normal subgroup of a (discrete) group $\, G\, $. For each such pair there is defined a sequence of abelian groups $\, \{ K_n^J ( N , G )\}\, ,\, n\geq 1\, $, which behave covariantly functorial for group morphisms $\, \phi :  ( N , G ) \rightarrow ( N' , G' )\, $, i.e. 
$\,\phi :  F \rightarrow F'\, $ is a homomorphism that restricts to $\, N \rightarrow N'\, $. The groups 
$\, K^J_n ( G ) = K^J_n ( G , G )\, $ are called {\it absolute} $K^J$-groups in contrast to the 
{\it relative} $K^J$-groups $\, K^J_n ( N , G )\, $ with $\, N\subset G\, $ a proper subgroup.
The regular theory on the other hand which is denoted $\, \{\, K_n ( N , G ) \}\, ,\, n\geq 1\, $ is only defined on a suitable subcategory, the category of \emph{ almost regular pairs} , which includes for example the full subcategory of normal subgroups of arbitrary perfect groups $\, G\, $. When both functors are defined there is a natural map $\, K_n ( N , G ) \rightarrow K_n^J ( N , G )\, $ which is injective. The name $K$-theory stems from the affinity with (higher) Algebraic $K$-theory (of Quillen, resp. Milnor in case of $K_2$), as well as topological and operator algebra $K$-theory. These together with group homology can be viewed as the main "neighbouring" theories. There is a natural map from $K^J$-theory to group homology (of degree zero) which is sort of an analogue to the Chern character in topological $K$-theory, resp. the corresponding map from Quillen's Algebraic $K$-theory to group homology which comes in the guise of the Hurewicz map from homotopy to homology groups by means of Quillens +-construction. Up to now the precise connection of the groups $\, \{ K_n ( R ) \}\, $ of a ring $\, R\, $ with Quillens higher algebraic $K$-groups is still unsettled in dimensions $\, n\geq 4\, $ whereas in dimension $\, n=2\, $ one has an isomorphism of the Regular $K_2$-group with Milnor's $\, K_2 ( R )\, $ which in turn coincides with Quillen's $\, K_2^Q ( R )\, $, and in dimension $\, n = 3\, $ there is a (surjective) map 
$\, K_3 ( R ) \twoheadrightarrow K^Q_3 ( R )\, $ the latter group which is equal to the third homology group of the Steinberg group $\, St ( R )\, $ from [Gersten] using the natural transformation
$$ B_n :\> K^J_n ( G )\>\longrightarrow\> H_n ( G ) = H_n ( G , \mathbb Z ) $$
mentioned above which is an isomorphism in dimension $\, n = 2\, $ and surjective in dimension 
$\, n=3\, $. It is more than likely that the induced transformation to Quillen $K$-theory for these special cases extends naturally to higher dimensions but this is a problem yet to be solved. The connection with topological $K$-theory is more obscure, first of all because up to now our theory is only defined for discrete structures and not for topological spaces per se. Still there is some congruence in spirit and many of the ideas employed in the construction are borrowed from corresponding concepts of topological $K$-theory, such as the concept of a mapping cone, suspension etc..
In addition one has an interesting connection of group $K^J$-theory 
with the topological $K$-homology of the classifying space $\, BG\, $ taking the form of a natural transformation
$$ C_{* + 2n} :\>  K^J_{* + 2n} ( G )\> \longrightarrow\> K^{top}_* ( BG ) \>  $$
for $\, * = 0 , 1\, $. 
A drawback of our definition of Algebraic $K$-theory might be seen in the fact that a priori the lowest (basic) dimension of the theory is $\, n = 2\, $. The higher $K$-groups are defined recurring to this dimension via an (algebraic) suspension operation, whereas the cases $\, n = 1\, $, and for rings, if one wishes $\, n \leq 0\, $, have to be defined separately. In case of $\, n = 1\, $ the definition is very simple (in general there is no distinction between the regular case and the broader $K^J$-groups in dimensions $\, n = 1 , 2\, $). It takes the form 
$$ K_1 ( N , G )\> =\> {N\over [\, N\, ,\, G\, ]} $$
where $\, [\, N\, ,\, G\, ]\, $ denotes the commutator subgroup of $\, N\, $ generated by commutators of the form $\, xyx^{-1}y^{-1}\, ,\, x\in N\, ,\, y\in G\, $, resp. 
$$ K_1 ( R )\> =\> {GL ( R )\over [\, GL ( R )\, ,\, GL ( R )\, ]}\> =\>  {GL ( R )\over E ( R )} $$
for rings which is just the usual definition of the first (relative) homology group of $\, ( N , G )\, $ 
(resp. first Algebraic $K$-group of $\, R\, $). Although many facets of the theory are still somewhat mysterious some fundamental results can be proved. One is the existence of long exact sequences associated with a group extension
$$ 0\>\longrightarrow N\>\longrightarrow G\>\longrightarrow\> H\>\longrightarrow\> 0\>  , $$
i.e. $\, N\, $ sits as a normal subgroup of $\, G\, $ with corresponding quotient $\, H = G / N\, $. For such an extension one may define boundary maps 
$\, {\delta }_n :\, K^J_{n+1} ( H ) \rightarrow K^J_n ( N , G )\, $ to obtain a long sequence of abelian groups and homomorphisms
$$ \cdots\longrightarrow\> K^J_{n+1} ( H )\>\buildrel {\delta }_n\over\longrightarrow\> 
K^J_n ( N , G )\>\longrightarrow\> K^J_n ( G )\>\longrightarrow\> K^J_n ( H )\> \buildrel 
{\delta }_{n-1}\over\longrightarrow\> \cdots $$ 
terminating in $\, K_1 ( H )\, $ such that the composition of any two consecutive homomorphisms is trivial. In fact it can be shown that the sequence is everywhere exact for general normal pairs with the possible exception of $\, K^J_3 ( N , G )\, $ and $\, K^J_3 ( G )\, $. One then has the notion of an {\it exact pair} 
$\, ( N , G )\, $ which property implies exactness at all places including $\, n=3\, $. Also for the regular 
$K$-groups there is an analogous long sequence which is everywhere exact in general. Another important property is {\it excision}, but this has a slightly different (weaker) meaning in our context than what is usually  understood by this term. Roughly it means that in higher dimensions ($\, n\geq 4\, $ is high enough) one has the following result: assume given two pairs $\, ( N , F )\, $ and $\, ( M , G )\, $ such that $\, F\subseteq G\, $ is a subgroup and $\, N = M\cap F\, $, then the functorial map
$\, K^J_n ( N , F ) \rightarrow K^J ( M , G )\, $ is an {\it injection}. In particular injective group homomorphisms induce injective maps on the level of "higher" absolute $K$-groups which certainly is a very strong result apt to arouse some suspicion. On the level of rings it also implies that on defining an 
"absolute" $K$-group of an ideal $\, I\subseteq R\, $ by considering its unitization $\, I^+\, $ adjoining the unit element of $\, R\, $ and taking the kernel of the maps 
$\, K_n ( I^+ ) \twoheadrightarrow K_n ( I^+ / I )\, $ one gets an injection 
$\, K_n ( I , I^+ ) \subseteq K_n ( I , R )\, $ for $\, n\geq 4\, $ (but not necessarily an isomorphism). 
There are other interesting results which may arouse suspicion among experts leading them to the conviction that the theory is altogether trivial  in higher dimensions ($\, n\geq 4\, $) as indeed it is for certain classes of groups containing both the classes of nilpotent and of finite groups. The author on the other hand is more optimistic and inclined to believe because of the existence of long exact sequences for the regular theory and by the belief that the functorial map $\, K_3 ( N , G )\rightarrow K_3 ( G )\, $ will not be injective in general, that there exist interesting classes of groups having nontrivial higher $K^J$-groups. The somewhat simplest such candidate would be an infinite twostep solvable group. If there is something good about the known vanishing results in higher dimensions it is that they add to the general computability of the theory and indicate that the higher $K$-groups are not "too large" and measure only "serious things", although built from a huge apparatus increasing exponentially in size with the dimension. From the finite group case one gets a corresponding result for finite rings $\, R\, $. Namely, in this case the group $\, E ( R )\, $, which itself is infinite, is the inductive limit of the sequence of finite groups $\, E_n ( R )\, $ generated by $n$-dimensional elementary matrices. Since the $K$-groups commute with inductive limits one finds that $\, K_n ( R ) = 0\, $ for any finite ring $\, R\, $ and 
$\, n\geq 4\, $. In particular this shows that our theory differs from Quillen's on considering the finite rings 
$\, {\mathbb Z}_m = \mathbb Z / m\mathbb Z\, $.
We now explain in more detail how the paper which consists of six sections and three appendices is organized. The first section stakes out the field giving some elementary constructions used throughout the paper. The definition of the relative $K$-groups in the basic dimension $\, n = 2\, $ is given along with some lemmas of mainly technical interest. The concept of suspension in $K^J$-theory is motivated
by looking at certain canonical free resolutions of the pair $\, ( N , G )\, $ but the explicit definitions of higher $K^J$-groups and suspension is postponed to later sections. Section 2 is devoted to the concept of regular suspension which takes the place of ordinary suspension in case of the regular theory.
Theorems 1 and 2 are motivating for the regular theory and exhibit the specific property of the regular groups, properly defined in section 6, which is their "smoothness"  or infinite regularity as opposed to $K^J$-theory which also measures things of finite (k-)regularity which do not reproduce themselves passing to finite commutator subgroups. Theorem 2 asserts that certain universal (and essentially unique) central extensions can be constructed on the iterated regular suspensions of an almost regular pair which then serves in later sections to show that the regular theory is well defined and exact. However the proof depends on some intricate technicalities and general results of the later sections so that in section 2 only a sketch of proof can be given for Theorem 2. Section 3 is an excursion to the topic of {\it strict splitting} which is supposed to be a property of a normal subgroup $\, J\, $ of a free group 
$\, U\, $ which  is reproduced by forming suspensions and entails trivial $K^J$-groups in all dimensions $\, n\geq 2\, $. This is of extreme importance since it allows to identify certain "conelike" $K$-trivial objects and to define {\it mapping cones} for morphisms of normal pairs. Section 4 is concerned with the above mentioned and overall important property of excision. The iterated suspensions and the higher $K^J$-groups are formally defined. Section 5 then gives the major part of exactness results for $K^J$-theory, i.e. the boundary maps and long sequence associated with a group extension are defined and investigated in detail to which instance this sequence is exact. Section 6 is the longest one and primarily about the regular theory, its main theorem being the long exact sequences in analogy with section 5. Since Regular $K$- and $K^J$-theory are so intimately interwoven the section also contains a lot of general results some which are used throughout the paper even in the earlier sections. Appendix A deals with the case of finite groups and proves the vanishing result in higher dimensions mentioned above. Appendix B is concerned with the construction of the natural transformation to group homology and Appendix C gives the construction of the mentioned natural transformation of the $K^J$-theory of 
$\, G\, $ to the topological $K$-homology of the classifying space $\, BG\, $. For $\, n=2\, $ this result has been applied in [Haag2] to prove rational injectivity of the analytic assembly map of Kasparov in dimension $\leq 2$. Finally let us say something about the future prospects of the theory. There are a lot of unsolved questions if not mysteries, hopefully provoking the curiosity of the interested reader, and some directions to which the theory may be extended. One such extension would be to the category of topological groups, another to monoids or categories. It would also be desirable to have a dual 
Co-$K$-theory behaving well in certain respects (existence of long exact sequences, products etc.) and indications are there that such a theory can be constructed, though maybe not in a unique or canonical way. It would also be helpful to compute some more nontrivial specific examples. Each such computation increases the computability of other cases by invoking long exact sequences, excision arguments and other general properties.
\par\bigskip\noindent
{\sc References}
\par\bigskip\noindent
[Atiyah] $\> $ M. Atiyah, {\it K-Theory}, Lecture notes by D. W. Anderson, 
\par\noindent 
$\quad $ Benjamin New York / Amsterdam (1967);
\par\smallskip\noindent
[Bass] $\> $ H. Bass, {\it Algebraic K-Theory}, Benjamin
\par\noindent
$\quad $  New York / Amsterdam (1968);
\par\smallskip\noindent
[Cuntz] $\> $ J. Cuntz, {\it Bivariant $K$-Theory and the Weyl algebra}, 
\par\noindent
$\quad $ Math. arXiv (2004);
\par\smallskip\noindent
[Cortinas-Thom] $\> $ G. Cortinas, A. Thom, {\it Bivariant Algebraic} 
\par\noindent
$\quad $ {\it $K$-Theory}, Math. arXiv (2007);
\par\smallskip\noindent
[Gersten] $\> $ S. M. Gersten, {\it $K_3$ of a ring is $H_3$ of the Steinberg group}, 
\par\noindent
$\quad $ Proc. Am. Math. Soc. {\bf 37} (1973);
\par\smallskip\noindent
[Grothendieck] $\> $ A. Borel, J.-P. Serre, {\it Le th\'eoreme de Riemann-Roch}, 
\par\noindent
$\quad $ Bull. Soc. Math. France {\bf 86} (1958);
\par\smallskip\noindent
[Haag1] $\> $ U. Haag, {\it Regular Algebraic $K$-Theory for Groups}, 
\par\noindent
$\quad $ 6 sections, 3 appendices, Preprint (2004--2010);
\par\smallskip\noindent
[Haag2] $\> $ U. Haag, {\it On rational injectivity of Kasparov's assembly} 
\par\noindent
$\quad $ {\it map in dimension $\leq 2$}, Math. arXiv (2012);
\par\smallskip\noindent
[Karoubi-Villamayor] $\> $ M. Karoubi, O. Villamayor, {\it $K$-theorie}  
\par\noindent
$\quad $ {\it algebrique et $K$-theorie topologique}, C. R. Acad. Sci. Paris 
\par\noindent
$\quad $ {\bf 269} (1969);
\par\smallskip\noindent
[Kasparov] $\> $ G. Kasparov, {\it The operator $K$-functor and extensions}  
\par\noindent
$\quad $ {\it of $C^*$-algebras}, Izv.Akad.Nauk SSSR Ser. Math. {\bf 44}:3 (1980);
\par\smallskip\noindent
[Milnor]  $\> $ J. Milnor, {\it Algebraic $K$-Theory and Quadratic Forms}, 
\par\noindent
$\quad $ Inventiones math. {\bf 9} (1970);
\par\smallskip\noindent
[Quillen] $\> $ D. Quillen, {\it Higher Algebraic $K$-Theory I},  Algebraic
\par\noindent
$\quad $ $K$-Theory I,  LNM {\bf 341} (1973);
\par\smallskip\noindent
[Schur] $\> $ J. Schur, {\it Ueber die Darstellungen der endlichen Gruppen}
\par\noindent
$\quad $ {\it durch gebrochene lineare Substitutionen}, J. Reine Angew. Math. 
\par\noindent
$\quad $  {\bf 127} (1904);
\par\smallskip\noindent
[Swan] $\> $ R. Swan, {\it Non-Abelian homological algebra and $K$-Theory}, 
\par\noindent
$\quad $ Proc. Symp. Pure Math., Am. Math. Soc. {\bf 17} (1970);
\par\noindent

\end{document}